\documentclass[12pt]{article}
\usepackage{a4wide}

\usepackage{amsmath,amssymb,amsthm}
\usepackage[mathscr]{eucal}
\usepackage{graphics}
\usepackage[hypertex]{hyperref}
\usepackage{epsfig}
\usepackage{psfrag}
\usepackage[usenames]{color}

\author{Ievgen~V.~Bondarenko, Rostyslav~V.~Kravchenko\thanks{The author was partially supported by NSF grants
DMS-0308985, DMS-0456185 and DMS-0600975}}
\title{\textbf{Graph-directed systems and self-similar measures on limit spaces\\ of self-similar groups}}

\newcommand{\xmo}{X^{-\omega}}
\newcommand{\xo}{X^\omega}
\newcommand{\xz}{X^{\mathbb{Z}}}
\newcommand{\Z}{\mathbb{Z}}
\newcommand{\Rs}{\hat{\mathbb{C}}}
\newcommand{\M}{\mathcal{M}}

\newcommand{\sol}[1][G]{\mathcal{S}_{#1}}
\newcommand{\limGs}[1][G]{\mathcal{X}_{#1}}
\newcommand{\lims}[1][G]{\mathscr{J}_{#1}}
\newcommand{\si}{\mathsf{s}}
\newcommand{\e}{\mathsf{e}}

\newcommand{\muls}{\mathsf{m}}
\newcommand{\mgr}{meas}

\newcommand{\tile}{\mathscr{T}}
\newcommand{\T}{\mathsf{T}}
\newcommand{\nucl}{\mathcal{N}}

\newcommand{\gr}{\Gamma}
\newcommand{\B}{\mathscr{B}}
\newcommand{\F}{\mathscr{F}}

\newtheorem{theor}{Theorem}
\newtheorem{prop}[theor]{Proposition}
\newtheorem{cor}[theor]{Corollary}
\newtheorem{lemma}{Lemma}
\theoremstyle{definition}
\newtheorem{defi}{Definition}
\newtheorem{example}{Example}
\newtheorem{remark}{Remark}

\begin{document}
\maketitle

\begin{abstract}
Let $G$ be a group and $\phi:H\rightarrow G$ be a contracting homomorphism from a subgroup $H<G$ of
finite index. V.~Nekrashevych \cite{self_sim_groups} associated with the pair $(G,\phi)$ the limit
dynamical system $(\lims,\si)$ and the limit $G$-space $\limGs$ together with the covering
$\cup_{g\in G}\tile\cdot g$ by the tile $\tile$. We develop the theory of self-similar measures
$\muls$ on these limit spaces. It is shown that $(\lims,\si,\muls)$ is conjugated to the one-sided
Bernoulli shift. Using sofic subshifts we prove that the tile $\tile$ has integer measure and we
give an algorithmic way to compute it. In addition we give an algorithm to find the measure of the
intersection of tiles $\tile\cap (\tile\cdot g)$ for $g\in G$. We present applications to the
invariant measures for the rational functions on the Riemann sphere and to the evaluation of the
Lebesgue measure of integral self-affine tiles.\\


\end{abstract}



\section{Introduction}

Let $G$ be a group and $\phi$ be a virtual endomorphism of $G$, which is a homomorphism from a
subgroup $H<G$ of finite index to $G$. Iterative construction involving $\phi$ (together with some
additional data) produces the so-called self-similar action $(G,X^{*})$ of the group $G$ on the
space $X^{*}$ of finite words over an alphabet $X$. And conversely, every self-similar action of
the group $G$ defines a virtual endomorphism of $G$, which almost completely describes the action.
A rich geometric theory is associated with the pair $(G,\phi)$ in \cite{self_sim_groups} through
the theory of self-similar groups. The goal of this paper is to introduce measure to this theory.

Self-similar group is a rather new notion in geometric group theory. Like the self-similar objects
in geometry (fractals) are too irregular to be described using the language of classical Euclidean
geometry, the self-similar groups possess properties not typical for the traditional group theory.
In particular, the class of self-similar groups contains infinite periodic finitely generated
groups, groups of intermediate growth, just-infinite groups, groups of finite width, etc. (see
\cite{self_sim_groups,fractal_gr_sets,branch,gri:MilnorProb} and references therein). At the same
time, it was discovered that self-similar groups appear naturally in many areas of mathematics, and
have applications to holomorphic dynamics, combinatorics, analysis on fractals, etc. An important
class of self-similar groups are contracting groups, which correspond to self-similar actions with
contracting virtual endomorphism. A virtual endomorphism $\phi$ is contracting if it asymptotically
contracts the length of group elements with respect to some generating set. The contracting
property makes many problems around the group effectively solvable.



V.~Nekrashevych in \cite{self_sim_groups} associated a limit dynamical system $(\lims,\si)$ with
every contracting self-similar action, where $\lims$ is a compact metrizable space and
$\si:\lims\rightarrow\lims$ is an expanding continuous map. The limit space $\lims$ can be defined
as the quotient of the space of left-infinite sequences $\xmo=\{\ldots x_2x_1 | x_i\in X\}$ by the
equivalence relation, which can be recovered from a finite directed labeled graph $\nucl$, called
the nucleus of the action. Another associated geometric object is the limit $G$-space $\limGs$,
which is a metrizable locally compact topological space with a proper co-compact (right) action of
$G$. The limit spaces $\lims$ and $\limGs$ depend up to homeomorphism only on the pair $(G,\phi)$.
However every self-similar action with the pair $(G,\phi)$ additionally produces a tile $\tile$ of
the limit $G$-space and a covering $\limGs=\cup_{g\in G}\tile\cdot g$ (not a tiling in general).

Limit spaces connect self-similar groups with the classical self-similar sets. The self-similar set
(the attractor) given by the system of contracting similarities $f_1,\ldots,f_n$ (iterated function
system) of a complete metric space is the unique compact set $T$ satisfying $T=\cup_{i=1}^n
f_i(T)$. Given a probability vector $p=(p_1,\ldots,p_n)$, Hutchinson
\cite{hutch:fract_self_similarity} showed the existence of a unique probability measure $\mu$
supported on $T$ satisfying
\begin{equation}\label{eqn_self_similar_measure}
\mu(A)=\sum_{i=1}^{n} p_i\mu(f_i^{-1}(A)), \mbox{ for any Borel set $A$},
\end{equation}
which is called the self-similar measure. Another way to introduce this measure is to consider the
natural coding map $\pi:\xmo\rightarrow T$ given by $\pi(\ldots x_2x_1)=\cap_{m\geq 1} f_{x_1}\circ
f_{x_2}\circ\ldots\circ f_{x_m}(T)$. Then the self-similar measure $\mu$ is the image of the
Bernoulli measure $\mu_p$ on $\xmo$ with weight $p$ (here $\mu_p(x_i)=p_i$) under the projection
$\pi$. Self-similar measures play an important role in the development of fractal geometry, and
have applications in harmonic analysis, conformal dynamics, algebraic number theory, etc. (see
\cite{edgar:int_prob_fr_mes,urb:meas_dim_confdyn, ifs_overl_self_sim_mes,str:self_harm,
bandt:self_sim_meas} and references therein).


The Bernoulli measure $\mu_p$ is one of the basic self-similar measures on the spaces of (left or
right) infinite sequences over the alphabet $X$. In Section~\ref{section_Sofic} we study the
Bernoulli measure of sofic subshifts and other sets given by a finite directed graph $\gr=(V,E)$,
whose edges are labeled by elements of $X$. Consider the set $F_v$ for $v\in V$ of all sequences
$\ldots x_2x_1$, which are read along left-infinite paths ending in the vertex $v$. It is proved in
Section~\ref{section_Sofic} that if the graph $\gr$ is right-resolving (i.e. for every vertex $v$
the outgoing edges at $v$ are labeled distinctly) then the sum $\mgr(\gr)=\sum_{v\in V} \mu_p(F_v)$
is integer, which does not depend on the probability vector $p$. It can be interpreted as follows:
almost every left-infinite sequence belongs to precisely $\mgr(\gr)$ sets $F_v$. The number
$\mgr(\gr)$ we call the measure number of the graph $\gr$. We propose an algorithmic method to
compute the measures $\mu_p(F_v)$ for any graph $\gr$, in particular its measure number.

The push-forward of the uniform Bernoulli measure on $\xmo$ provides the self-similar measure
$\muls$ on the limit space $\lims$. The $G$-invariant measure $\mu$ on the limit $G$-space $\limGs$
is defined in a similar way. The measure $\mu$ restricted to the tile $\tile$ satisfies the
self-similarity equation (\ref{eqn_self_similar_measure}), so it is also self-similar. It is proved
in Section~\ref{section_Measures} that the measure $\mu(\tile)$ is equal to the measure number of
the nucleus $\nucl$. In particular it is integer, the fact which generalizes corresponding result
for integral self-affine tiles \cite{lag_wang:int_self_affine_I}. In addition we give an algorithm
to find the measure of intersection of tiles $\tile\cap (\tile\cdot g)$ for $g\in G$. Then the
covering $\limGs=\cup_{g\in G}\tile\cdot g$ is a perfect multiple covering of multiplicity
$\mu(\tile)$, i.e. every point of $\limGs$ belongs to at least $\mu(\tile)$ tiles and almost every
point belongs to precisely $\mu(\tile)$ tiles. This is used to prove that the measures $\muls$ and
$\mu$ depend not on the specific self-similar action of $G$, but only on the pair $(G,\phi)$ as the
limit spaces themselves. Using a criterion from \cite{hoffman:Bernoulli} we show that the limit
dynamical system $(\lims,\si,\muls)$ is conjugated to the one-sided Bernoulli shift.

This work is partially motivated by applications presented in Section~\ref{section_Applictions}. If
$G$ is a torsion-free nilpotent group with a contractive surjective virtual endomorphism $\phi$ and
a faithful self-similar action, then the measure $\mu$ on $\limGs$ can be considered as a Haar
measure on the respective nilpotent Lie group, Malcev's completion of $G$. In the case of
self-similar actions of the free abelian group $\mathbb{Z}^n$ the limit $G$-space
$\limGs[\mathbb{Z}^n]$ is $\mathbb{R}^n$ and the tile $\tile$ is an integral self-affine tile,
which are intensively studied for the last two decades (see \cite{lag_wang:self_affine,
vince:digit, lag_wang:int_self_affine_I,lag_wang:int_self_affine_II,he_hui:self_affine}). In this
case the measure $\mu$ is the Lebesgue measure on $\mathbb{R}^n$. One can apply the methods
developed in Section~\ref{section_Sofic} to give an algorithmic way to find the Lebesgue measure of
an integral self-affine tile, providing answer to the question in \cite{lag_wang:int_self_affine_I}
(initially solved in \cite{gab_yu:natur_til} without self-similar groups). In addition we have an
algorithm to find the Lebesgue measure of the intersection of tiles $\tile\cap (\tile+a)$ for
$a\in\mathbb{Z}^n$ studied in \cite{gab_yu:natur_til,int_two_tiles}. Finally, the theory of
iterated monodromy groups developed in \cite{self_sim_groups} implies that if $G$ is the iterated
monodromy group of a sub-hyperbolic rational function $f$ then the measure $\muls$ on the limit
space $\lims$ is the unique $f$-invariant probability measure $\mu_f$ on the Julia set of $f$
studied in \cite{invmeas_rat,lub:entr_ration,rational_bernoulli}.

\section{Bernoulli measure of sofic subshifts}\label{section_Sofic}

Let $X$ be a finite set with discrete topology. Denote by $X^{*}=\{x_1x_2\ldots x_n | x_i\in X,
n\geq 0\}$ the set of all finite words over $X$ (including the empty word denoted $\emptyset$). Let
$X^{\omega}$ be the set of all right-infinite sequences (words) $x_1x_2\ldots$, $x_i\in X$. Let
$X^{-\omega}$ be the set of all left-infinite sequences (words) $\ldots x_2x_1$, $x_i\in X$. We put
the product topology on these sets. The length of a word $v=x_1x_2\ldots x_n$ is denoted by $|v| =
n$.

The \textit{shift on the space $\xo$} is the map $\sigma:\xo\rightarrow\xo$, which deletes the
first letter of a word, i.e. $\sigma(x_1x_2x_3\ldots) = x_2x_3\ldots$. The \textit{shift on the
space $\xmo$} is the map also denoted by $\sigma$, which deletes the last letter of a word, i.e.
$\sigma(\ldots x_3x_2x_1) = \ldots x_3x_2$. The shifts are continuous $|X|$\,-\,to\,-$1$ maps. The
branches $\sigma_x$ for $x\in X$ of the inverse of $\sigma$ are defined by
$\sigma_x(x_1x_2\ldots)=xx_1x_2\ldots$ and $\sigma_x(\ldots x_2x_1)=\ldots x_2x_1x$.

Let $p=(p_x)_{x\in X}$ be a probability vector (fixed for the rest of the section) and let $\mu_p$
be the \textit{Bernoulli measure} on $\xo$ with weight $p$, i.e. this measure is defined on
cylindrical sets by
\[
\mu_p(x_1x_2\ldots x_n\xo)=p_{x_1}p_{x_2}\ldots p_{x_n}.
\]
The measure on $\xmo$ is defined in the same way. We always suppose that $p_x>0$ for all $x\in X$
(otherwise we can pass to a smaller alphabet $X$). In case $p_x=\frac{1}{|X|}$ for all $x\in X$,
the measure $\mu_p$ is the \textit{uniform Bernoulli measure} denoted $\mu_u$. The dynamical system
$(\xo,\sigma,\mu_u)$ is called the \textit{one-sided Bernoulli $|X|$-shift}. The measure $\mu_p$ is
the unique regular Borel probability measure on $\xo$ that satisfies the self-similarity condition:
\[
\mu_p(A)=\sum_{x\in X} p_x\mu_p(\sigma_x^{-1}(A))
\]
for any Borel set $A\subset\xo$.


In this section all considered graphs are directed and labeled with $X$ as the set of labels. Let
$\gr=(V,E)$ be such a graph and take a vertex $v\in V$. We say that a sequence $x_1x_2\ldots\in\xo$
(a word $x_1x_2\ldots x_n\in X^{*}$) \textit{starts in the vertex $v$} if there exists a
right-infinite path $e_1e_2\ldots$ (finite path $e_1e_2\ldots e_n$) in $\gr$, which starts in $v$
and is labeled by $x_1x_2\ldots$ (respectively $x_1x_2\ldots x_n$). Similarly, we say that a
sequence $\ldots x_2x_1\in\xmo$ (a word $x_n\ldots x_2x_1\in X^{*}$) \textit{ends in the vertex
$v$} if there exists a left-infinite path $\ldots e_2e_1$ (finite path $e_n\ldots e_2e_1$) in
$\gr$, which ends in $v$ and is labeled by $\ldots x_2x_1$ (respectively $x_n\ldots x_2x_1$). For
every $w\in X^*\cup\xmo$ denote by $V_\gr(w)=V(w)\subset V$ the set of all vertices $v\in V$ such
that the sequence $w$ ends in $v$. Observe that $V(w'w)\subseteq V(w)$ for arbitrary word $w'$ and
finite word~$w$.

For every vertex $v\in V$ denote by $B_v$ the set of all right-infinite sequences that start in
$v$, and denote by $F_v$ the set of all left-infinite sequences that end in $v$. The sets $B_v$ and
$F_v$ are closed correspondingly in $\xo$ and $\xmo$, thus compact and measurable. The sets
$\B=\cup_{v\in V}B_v$ and $\F=\cup_{v\in V}F_v$ are the one-sided (respectively, right and left)
\textit{sofic subshifts} associated with the graph $\gr$. The sets $F_v, v\in V$, satisfy the
recursion
\[
F_v=\bigcup_{u\, \stackrel{x}{\rightarrow} \,v} \sigma_x(F_u)
\]
(here the union is taken over all edges which end in $v$). Hence, associating the map $\sigma_x$
with every edge $e$ of the graph $\gr$ labeled by $x$, the collection of sets $\{F_v, v\in V\}$ can
be seen as the graph-directed iterated function system on the sofic subshift $\F$ with the
underlying graph $\gr$ (see \cite{fractal_gr_sets}). All the maps $\sigma_x^{-1}$  are restrictions
of the shift $\sigma$, and thus $\{F_vx, v\in V\}$ is the Markov partition of the dynamical system
$(\F,\sigma)$. Similarly, the collection of sets $\{B_v, v\in V\}$ can be seen as the
graph-directed iterated function system on the sofic subshift~$\B$.

A labeled graph $\gr=(V,E)$ is called \textit{right-resolving} (Shannon graph in some terminology)
if for every vertex $v\in V$ the edges starting at $v$ have different labels. Every sofic subshift
can be given by a right-resolving graph (see Theorem~3.3.2 in \cite{IntroSymbDyn}). A
right-resolving graph is called \textit{strictly right-resolving} if every vertex $v\in V$ has an
outgoing edge labeled by $x$ for every $x\in X$.

%
%

For a labeled graph $\gr=(V,E)$ we use the following notations:
\[
\vec\mu_p(\B)=(\mu_p(B_v))_{v\in V}, \quad  \vec\mu_p(\F)=(\mu_p(F_v))_{v\in V}, \quad \mbox{ and }
\quad \mu_p(\gr)=\sum_{v\in V} \mu_p(F_v).
\]
Next we study the properties of these quantities and describe an algorithmic way to find them.
First we discuss the problem for right-resolving graphs, and then reduce the general case to the
right-resolving one.

\begin{theor}\label{theor_measure_graph_integer}
Let $\gr=(V,E)$ be a finite right-resolving graph. Then
\[
\mu_p(\gr)=\min_{w\in\xmo} |V(w)|=\min_{w\in X^*} |V(w)|.
\]
In particular, the measure $\mu_p(\gr)$ is integer.
\end{theor}
\begin{proof}
Let $w=\ldots x_2x_1\in\xmo$ and denote $w_n=x_n\ldots x_2x_1$ for $n\geq 1$. Observe, that
$V(w)\subseteq V(w_n)$ and $V(w_n)\subseteq V(w_m)$ for $n\geq m$.

Take a vertex $v\in \cap_{n\geq 1} V(w_n)$. Let $P_n$ be the set of all paths $e_n\ldots e_2e_1$
labeled by $w_n$ and ending in $v$. The set $P_n$ is a finite non-empty set for every $n$, and
$e_{n-1}\ldots e_2e_1\in P_{n-1}$ for every path $e_n\ldots e_2e_1\in P_n$. Since the inverse limit
of a sequence of finite non-empty sets is non-empty, there exists a left-infinite path $\ldots
e_2e_1$ labeled by $w$ and ending in $v$. Then $v\in V(w)$ and we get
\begin{equation}\label{eqn_th1_V(w) eq cap V(w_n)}
V(w)=\bigcap_{n\geq 1} V(w_n).
\end{equation}
From this follows that $|V(w)|=\min\limits_{n\geq 1} |V(w_n)|$ for $w\in\xmo$. Then
\[
\min_{w\in\xmo} |V(w)|=\min_{w\in\xmo} \min_{n\geq 1} |V(w_n)|=\min_{w\in X^*} |V(w)|
\]
and the second equality of the theorem is proved.

Define the integer $k=\min_{w\in\xmo} |V(w)|$ and consider the set
$\mathcal{O}=\mathcal{O}(\gr)\subseteq\xmo$ of all sequences $w\in\xmo$ such that $|V(w)|=k$.
Define $\mathcal{O}^*$ as the set of finite words that satisfy the same condition. We have the
following lemma.
\begin{lemma}\label{lemma_set V open dense}
The set $\mathcal{O}$ is open and dense in $\xmo$, and $\mu_p(\mathcal{O})=1$. For each
$w\in\mathcal{O}$ there is a beginning of $w$ that belongs to $\mathcal{O}^*$. Equivalently,
$\mathcal{O}=\cup_{w\in\mathcal{O}^*}\xmo w$.
\end{lemma}
\begin{proof}
If $w\in\mathcal{O}$ then $k=|V(w)|=\min_{n\geq 1} |V(w_n)|$ and there exists $N\geq 1$ such that
$|V(w_N)|=k$. Then $k\leq |V(\omega w_N)|\leq |V(w_N)|=k$ for all $\omega\in\xmo$. Hence
$w_N\in\mathcal{O}^*$, $\xmo w_N\subseteq\mathcal{O}$ and so $\mathcal{O}$ is open, and thus
measurable.

Let $u\in X^*$ be such that $|V(u)|=k$. Let us show that if $w\in\xmo$ contains the subword $u$
then $w\in\mathcal{O}$. If $w=w'u$ then $k\leq |V(w'u)|\leq |V(u)|=k$ and $w\in\mathcal{O}$.
Observe that $V(ux)$ is the set of those vertices $v\in V$ for which there exists an edge labeled
by $x$ which starts in some vertex of $V(u)$ and ends in $v$. Since the graph $\gr$ is
right-resolving there is no more than one such an edge for each vertex of $V(u)$, and thus
$|V(ux)|\leq |V(u)|$. It implies that if $w=w'uu'$ than $k\leq |V(w)|\leq |V(u)|=k$ and thus
$|V(w)|=k$, so $w\in\mathcal{O}$. The Bernoulli measure of the set of all words $w'uu', u'\in X^*,
w'\in\xmo$, is equal to $1$. Thus $\mu_p(\mathcal{O})=1$. It follows also that $\mathcal{O}$ is
dense in $\xmo$.
\end{proof}

By construction of the set $\mathcal{O}$ for every $w\in\mathcal{O}$ there exist exactly $k$
vertices $v$ such that $w\in F_v$. Let $\chi_{F_v}$ be the characteristic function of the set
$F_v$. Then $\sum_{v\in V}\chi_{F_v}=k$ almost everywhere. Integrating we get
$\mu_p(\gr)=\sum\limits_{v\in V} \mu_p(F_v)=k$.
\end{proof}

%

\begin{remark} The theorem holds not only for Bernoulli measures. The only property that was used is that
for every word $u\in X^{*}$ the set of all words, which contain $u$ as a subword, has measure $1$.
It should be pointed out that the number $\mu(\gr)$ is independent on the chosen measure $\mu$,
while the measures $\mu(F_v)$, $\mu(B_v)$ and $\sum_{v\in V} \mu(B_v)$ depend on $\mu$.
\end{remark}

\begin{defi}
The number $\mgr(\gr)=\mu_p(\gr)$ is called the \textit{measure number} of the graph~$\gr$.
\end{defi}

Theorem~\ref{theor_measure_graph_integer} shows that almost every sequence $w\in\xmo$ ends in
precisely $\mgr(\gr)$ vertices of $\gr$.

We will use the following proposition in the next sections.
\begin{prop}\label{prop cup B_v eq cup F_v}
Let $\gr=(V,E)$ be a finite labeled graph. Then $\B=\xo$ if and only if $\F=\xmo$. In particular,
$\F=\xmo$ for a finite strictly right-resolving graph $\gr$.
\end{prop}
\begin{proof}
Since the inverse limit of nonempty finite sets is nonempty, $\B=\xo$ ($\F=\xmo$) is equivalent to
the fact that every finite word $v\in{X^*}$ labels some path in $\gr$.
\end{proof}

The matrix $A=(a_{vu})_{v,u\in V}$, where $a_{vu}$ is equal to the number of edges from $v$ to $u$,
is the \textit{adjacency matrix} of the graph $\gr$. For the probability vector $p=(p_x)_{x\in X}$
define the matrix
\[
T_p=(t_{vu})_{v,u\in V}, \mbox{ where } t_{vu}=\sum\limits_{v\stackrel{x}\rightarrow u}p_x
\]
(the sum is taken over all edges from $v$ to $u$). The matrices $A$ and $T_p$ are irreducible if
and only if the graph $\gr$ is strongly connected. If the graph $\gr$ is right-resolving, then the
matrix $T_p$ is the \textit{transition matrix} of the random walk on the weighted directed graph
$\gr$, where each edge labeled by $x$ has weight $p_x$. In this case the row sums of the matrix $A$
are $\leq |X|$, and the row sums of the matrix $T_p$ are $\leq 1$, hence the spectral radius of $A$
is $\leq |X|$, and the spectral radius of $T_p$ is $\leq 1$. If the graph $\gr$ is strictly
right-resolving, then the transition matrix $T_p$ is right stochastic.

\begin{prop}\label{prop_eigenvect_transition_matr}
Let $\gr=(V,E)$ be a finite right-resolving graph with the transition matrix~$T_p$. If the vector
$\vec \mu_p(\B)$ is nonzero then it is the right eigenvector of $T_p$ for the eigenvalue~$1$. If
the vector $\vec \mu_p(\F)$ is nonzero then it is the left eigenvector of $T_p$ for the
eigenvalue~$1$.
\end{prop}
\begin{proof}
By construction, for every vertex $v\in V$ we have
\[
B_v=\bigsqcup_{v\stackrel{x}\rightarrow u} xB_u
\]
(here the union is disjoint because the graph $\gr$ is right-resolving). It implies
\[
\mu_p(B_v)=\sum_{v\stackrel{x}\rightarrow u} p_x\mu_p(B_{u})=\sum_{u\in V} t_{vu}\mu_p(B_{u}).
\]
Thus the nonzero vector $\vec \mu_p(\B)$ is the right eigenvector of $T_p$ for the eigenvalue~$1$.

Similarly, 
\[
F_v=\bigcup_{u\stackrel{x}\rightarrow v} F_{u}x,\ v\in V,
\]
and, since the graph $\gr$ is right-resolving, that implies
\begin{equation*}
\mu_p(F_v)\leq\sum_{u\stackrel{x}\rightarrow v} p_x\mu_p(F_{u})=\sum_{u\in V}
t_{uv}\mu_p(F_{u})\quad \Rightarrow\quad \vec \mu_p(\F)\leq\vec \mu_p(\F) T_p.
\end{equation*}

The standard arguments based on the theory of nonnegative matrices (see for example \cite[proof of
Theorem 4.5]{he_hui:self_affine}, \cite[p. 197]{gab_yu:natur_til}, \cite{handbook:rothblum}) end
the proof.
\end{proof}

\begin{cor}\label{cor_meas_sum_meas}
Let $\gr=(V,E)$ be a finite right-resolving graph. Let $\{\gr_i\}$ be the set of all strongly
connected components of $\gr$, which are strictly right-resolving graphs. Then
$\mgr(\gr)=\sum_i\mgr(\gr_i)$.
\end{cor}

In particular, if a finite strictly right-resolving graph $\gr$ contains a vertex $v_0$ such that
for each vertex $v$ there is a path in $\gr$ from $v$ to $v_0$ and for each $x\in X$ there is an
edge from $v_0$ to $v_0$ labeled by $x$ (the open set condition for graphs), then $\mgr(\gr)=1$.

\begin{cor}
Let $\gr=(V,E)$ be a finite right-resolving graph with the adjacency matrix~$A$. For the uniform
Bernoulli measure $\mu_u$ the nonzero vectors $\vec \mu_u(\B)$ and $\vec \mu_u(\F)$ are
respectively the right and left eigenvectors of $A$ for the eigenvalue $|X|$.
\end{cor}

\begin{figure}
\begin{center}
\psfrag{0}{\small{$0$}} \psfrag{1}{\small{$1$}} \psfrag{a}{\small{$a$}} \psfrag{b}{\small{$b$}}
\psfrag{e}{\small{$a_1$}} \psfrag{f}{\small{$b_1$}} \psfrag{c}{\small{$a_2$}}
\psfrag{d}{\small{$b_2$}} \epsfig{file=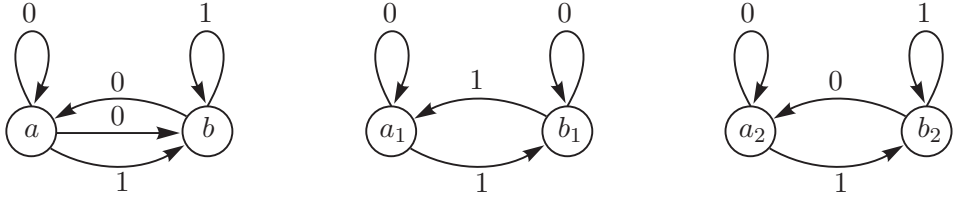,height=70pt}\caption{The graphs $\gr$ (on the
left), $\gr_1$ (in the center), and $\gr_2$ (on the right)} \label{fig_Example_3_Graphs}
\end{center}
\end{figure}

The next example shows that the right-resolving property in
Proposition~\ref{prop_eigenvect_transition_matr} and Theorem~\ref{theor_measure_graph_integer}
cannot be dropped even in the case of the uniform Bernoulli measure. The second example presents
two right-resolving graphs with the same transition matrices, but with different measure numbers.
\begin{example}
For the graph $\gr$ shown in Figure~\ref{fig_Example_3_Graphs} we get $F_a=\xmo 0$, $F_b=\xmo$, the
vector $\vec\mu_u(\F)=(1/2,1)$ is not an eigenvector of the adjacency matrix, and $\mu_u(\gr)=3/2$.
\end{example}

\begin{example}
Consider two graphs $\gr_1$ and $\gr_2$ shown in Figure~\ref{fig_Example_3_Graphs}. They have the
same transition matrix (for the uniform Bernoulli measure), however $F_{a_1}=F_{b_1}=\xmo$ and
$F_{a_2}=\xmo 0$, $F_{b_2}=\xmo 1$. Hence, $\mu_u(\gr_1)=2$ and $\mu_u(\gr_2)=1$.
\end{example}

Although Theorem~\ref{theor_measure_graph_integer} gives a usefull characterization of the number
$\mu_p(\gr)$, it does not present an algorithmic way to find it. It follows from
Proposition~\ref{prop_eigenvect_transition_matr} that the problem of finding $\vec\mu_p(\F)$ and
$\mu_p(\gr)$ reduces to the strongly connected components which are strictly right-resolving graphs
(for all other vertices $\mu_p(F_v)=0$). Notice that if $\gr$ is a strongly connected strictly
right-resolving graph, then the vector $\frac{1}{\mu_p(\gr)}\vec \mu_p(\F)$ is the unique
stationary probability distribution of the stochastic matrix $T_p$.


At the same time Proposition~\ref{prop_eigenvect_transition_matr} implies the algorithm to find the
vector $\vec\mu_p(\B)$ for a right-resolving graph. Indeed, a left eigenvector of $T_p$ for the
eigenvalue $1$ is uniquely defined if we know its entries $\mu_p(B_v)$ for vertices $v$ in the
strongly connected components of $\gr$ without outgoing edges. For every such a component $\gr'$,
we have $B_v=\xo$ and $\mu_p(B_v)=1$ for every vertex $v\in\gr'$ if the component $\gr'$ is a
strictly right-resolving graph, and $\mu_p(B_v)=0$ otherwise. In particular, if the matrix $T_p$ is
rational then the values $\mu_p(B_v)$ are rational.


%
%
%

Given a finite labeled graph $\gr=(V,E)$ the problems of finding the measures $\mu_p(B_v)$ and
$\mu_p(F_v)$ are equivalent, and can be reduced to right-resolving graphs. The problem is that the
system
\[
B_v=\bigcup_{v\stackrel{x}\rightarrow u} xB_{u}=\bigsqcup_{x\in
X}x\left(\bigcup_{v\stackrel{x}\rightarrow u} B_{u}\right),\ v\in V
\]
is not self-similar in the sense that the above expression involves not only the sets $B_v$ but
also their finite unions. Introducing additional terms to this recursion corresponding to these
unions we get a system with a right-resolving graph. This procedure is similar to the construction
described in the proof of Theorem~3.3.2 in \cite{IntroSymbDyn}.

\begin{prop}
For every finite graph $\gr=(V,E)$ one can construct a finite right-resolving graph $\gr'=(V',E')$
with the property that for every $v\in V$ there exists $v'\in V'$ such that $B_v=B_{v'}$.
\end{prop}
%

Similarly, one can find the measures of subshifts $\B$ and $\F$ by introducing new vertices
corresponding to $\cup_{v\in V} B_v$ and $\cup_{v\in V} F_v$.


\begin{cor}\label{cor_meas_gr_rational}
Let $\gr=(V,E)$ be a finite labeled graph. For the uniform Bernoulli measure $\mu_u$ all measures
$\mu_u(B_v)$, $\mu_u(F_v)$, $\mu_u(\gr)$ are rational.
\end{cor}

Consider the question how to find the measure of the intersection $B_v\cap B_u$ for $v,u\in V$.
Construct a new graph $\mathscr{G}$ with the set of vertices $V\times V$ and put an edge from
$(v,u)$ to $(v',u')$ labeled by $x\in X$ for every edges $v\stackrel{x}\rightarrow v'$ and
$u\stackrel{x}\rightarrow u'$ in the graph $\gr$ (label products of graphs by Definition~3.4.8 in
\cite{IntroSymbDyn}). It is easy to see that then $B_{(v,u)}=B_v\cap B_u$ (see Proposition~3.4.10
in \cite{IntroSymbDyn}).

The described methods work in a more general case when labels are words over $X$. Let $\gr=(V,E)$
be a directed graph labeled by labels from $X^{*}\setminus\{\emptyset\}$. Define the sets $B_v$ and
$F_v$ for $v\in V$ as above. Similarly, the collection of sets $\{B_v, v\in V\}$ can be seen as the
graph-directed function system with the underlying graph $\gr$, where the map $\sigma_{x_1\ldots
x_n}=\sigma_{x_n}\circ\ldots\circ\sigma_{x_1}$ is associated with every edge labeled by $x_1\ldots
x_n$. To find the measures of the sets $B_v$ and $F_v$, we can replace each edge labeled by a word
$x_1x_2\ldots x_k$ by a sequence of edges labeled by letters $x_1,x_2,\ldots,x_k$, and apply the
above method to this new graph. Also we can directly use
Proposition~\ref{prop_eigenvect_transition_matr}, where such graph $\gr$ is called
\textit{right-resolving} if for every vertex $v\in V$ the word-label of an edge starting at $v$ is
not a prefix of the word-label of another edge starting at $v$ (i.e. in some sense the associated
automaton is deterministic). The transition matrix $T_p$ is defined similarly
\[
T_p=(t_{vu})_{v,u\in V}, \mbox{ where } t_{vu}=\sum\limits_{v\xrightarrow{x_1...x_k}
u}p_{x_1}p_{x_2}\ldots p_{x_k}.
\]


\section{Self-similar actions and its limit spaces}\label{section_Self-similar_groups}
We review in this section the basic definitions and theorems concerning self-similar groups. For a
more detailed account and for the references, see~\cite{self_sim_groups}.

\textbf{Self-similar actions.} A faithful action of a group $G$ on the set $X^*$ is called
\emph{self-similar} if for every $g\in G$ and $x\in X$ there exist $h\in G$ and $y\in X$ such that
\[
g(xw)=yh(w)
\]
for all $w\in X^*$. The element $h$ is called the \textit{restriction} of $g$ on $x$ and denoted
$h=g|_x$. Inductively one defines the restriction $g|_{x_1x_2\ldots x_n}=g|_{x_1}|_{x_2}\ldots
|_{x_n}$ for every word $x_1\ldots x_n\in X^n$. Notice that $(g\cdot h)|_v=g|_{h(v)}\cdot h|_v$ (we
are using left actions).


\textbf{Virtual endomorphisms.} The study of the self-similar actions of a group is in some sense
the study of the \textit{virtual endomorphisms} of this group, which are homomorphisms from a
subgroup of finite index to the group. There is a general way to construct a self-similar
representation of a group with a given associated virtual endomorphism. Let $\phi: H\rightarrow G$
be a virtual endomorphism of the group $G$, where $H<G$ is a subgroup of index $d$. Let us choose a
left coset transversal $T=\{g_0,g_1,\ldots,g_{d-1}\}$ for the subgroup $H$, and a sequence
$C=\{h_0,h_1,\ldots,h_{d-1}\}$ of elements of $G$ called a \emph{cocycle}. The self-similar action
$(G,X^{*})$ with the alphabet $X=\{x_0,x_1,\ldots,x_{d-1}\}$ defined by the triple $(\phi,T,C)$ is
given by
\[
g(x_i)=x_j,\qquad g|_{x_i}=h_j^{-1}\phi(g_j^{-1}gg_i)h_i,
\]
where $j$ is such that $g_j^{-1}gg_i\in H$ (such $j$ is unique). The action may be not faithful,
the kernel can be described using Proposition~2.7.5 in \cite{self_sim_groups}.

Conversely, every self-similar action can be obtained in this way. Let $(G,X^{*})$ be a
self-similar action and take a letter $x\in X$. The stabilizer $St_G(x)$ of the letter $x$ in the
group $G$ is a subgroup of index $\leq |X|$ in $G$. Then the map $\phi_x: g\mapsto g|_x$ is a
homomorphism from $St_G(x)$ to $G$ called \emph{the virtual endomorphism associated to the
self-similar action}. Choose $T=\{g_y : y\in X\}$ and $C=\{h_y : y\in X\}$ such that $g_y (x)=y$
and $h_y=\left(g_y|_x\right)^{-1}$. Then $T$ is a coset transversal for the subgroup $St_G(x)$ and
the self-similar action $(G,X^*)$ is defined by the triple $(\phi_x,T,C)$. Different self-similar
actions of the group $G$ with the same associated virtual endomorphism are conjugated by
Proposition~2.3.4 in \cite{self_sim_groups}.



\textbf{Contracting self-similar actions.} An important class of self-similar actions are
contracting actions. A self-similar action of a group $G$ is called \emph{contracting} if there
exists a finite set $\nucl$ such that for every $g\in G$ there exists $k\in\mathbb{N}$ such that
$g|_v\in\nucl$ for all words $v\in X^*$ of length $\geq k$. The smallest set $\nucl$ with this
property is called the \emph{nucleus} of the self-similar action. The nucleus itself is
self-similar in the sense that $g|_v\in\nucl$ for every $g\in\nucl$ and $v\in X^{*}$. It can be
represented by the \textit{Moore diagram}, which is the directed labeled graph with the set of
vertices $\nucl$, where there is an edge from $g$ to $g|_x$ labeled $(x,g(x))$ for every $x\in X$
and $g\in\nucl$. We identify the nucleus with its Moore diagram, also denoted by $\nucl$. The
contracting property of the action depends only on the virtual endomorphism but not on the chosen
coset transversal and cocycle (see Corollary~2.11.7 in \cite{self_sim_groups}). Notice that every
contracting self-similar group is countable.

Self-similar groups are related to self-similar sets through the notion of limit spaces.

\textbf{Limit $G$-spaces.} Let us fix a contracting self-similar action $(G,X^{*})$. Consider the
space $\xmo\times G$ of all sequences \ $\ldots x_2x_1\cdot g$, $x_i\in X$ and $g\in G$, with the
product topology of discrete sets $X$ and $G$. Two elements $\ldots x_2x_1\cdot g$ and $\ldots
y_2y_1\cdot h$ of $\xmo\times G$ are called \textit{asymptotically equivalent} if there exist a
finite set $K\subset G$ and  a sequence $g_n\in K, n\geq 1,$ such that
\begin{displaymath}
g_n(x_n x_{n-1}\ldots x_1)=y_n y_{n-1}\ldots y_1\qquad \mbox{ and }\qquad g_n|_{x_n x_{n-1}\ldots
x_1}\cdot g=h
\end{displaymath}
for every $n\geq 1$. This equivalence relation can be recovered from the nucleus $\nucl$ of the
action (Proposition~3.2.6 in \cite{self_sim_groups}).

\begin{prop}\label{prop_asympt_Moore_diagr}
Two elements $\ldots x_2x_1\cdot g$ and $\ldots y_2y_1\cdot h$ of $\xmo\times G$ are asymptotically
equivalent if and only if there exists a left-infinite path $\ldots e_2e_1$ in the nucleus $\nucl$
ending in the vertex $hg^{-1}$ such that the edge $e_i$ is labeled by $(x_i,y_i)$.
\end{prop}

The quotient of the set $\xmo\times G$ by the asymptotic equivalence relation is called the
\textit{limit $G$-space} of the action and denoted $\limGs[(G,X^*)]$. The group $G$ naturally acts
on the space $\limGs[(G,X^*)]$ by multiplication from the right.

The map $\tau_x$ defined by the formula
\[
\tau_x(\ldots x_2 x_1\cdot g)=\ldots x_2 x_1 g(x)\cdot g|_x
\]
is a well-defined continuous map on the limit $G$-space $\limGs$ for every $x\in X$, which is not a
homeomorphism in general. Inductively one defines $\tau_{x_1x_2\ldots
x_n}=\tau_{x_n}\circ\tau_{x_{n-1}}\circ\ldots\circ\tau_{x_1}$.


The image of $\xmo\times 1$ in $\limGs$ is called the \textit{(digit) tile} $\tile$ of the action.
The image of $\xmo v\times 1$ for $v\in X^n$ is called the \textit{tile} $\tile_v$, equivalently
$\tile_v=\tau_v(\tile)$. It follows directly from definition that
\[
\limGs=\bigcup_{g\in G} \tile\cdot g \ \ \ \ \mbox{ and } \ \ \ \ \tile=\bigcup_{v\in X^n} \tile_v.
\]
Two tiles $\tile\cdot g$ and $\tile\cdot h$ intersect if and only if $gh^{-1}\in\nucl$. A
contracting action $(G,X^*)$ satisfies the \emph{open set condition} if for any element $g$ of the
nucleus $\nucl$ there exists a word $v\in X^*$ such that $g|_v=1$, i.e. in the nucleus $\nucl$
there is a path from any vertex to the trivial state. If the action satisfies the open set
condition then the tile $\tile$ is the closure of its interior, and any two different tiles have
disjoint interiors; otherwise every tile $\tile\cdot g$ is covered by the other tiles (see
Proposition~3.3.7 in \cite{self_sim_groups}).


The tile $\tile$ and the partition of $\limGs$ on tiles $\tile\cdot g$ depend on the specific
self-similar action of the group $G$. However, up to homeomorphism the limit $G$-space
$\limGs[(G,X^*)]$ is uniquely defined by the associated virtual endomorphism $\phi$ of the group,
hence we denote it by $\limGs[G](\phi)$ (or $\limGs$ for short).

\begin{theor}\label{theor_limsX hom limsY}
Let $\phi: H\rightarrow G$ be a virtual endomorphism of the group $G$. Let $(G,X^*)$ and $(G,Y^*)$
be the contracting self-similar actions defined respectively by the triples $(\phi,T,C)$ and
$(\phi,T',C')$. Then $\limGs[(G,X^*)]$ and $\limGs[(G,Y^*)]$ are homeomorphic and the homeomorphism
is the map $\alpha : \limGs[(G,X^*)] \longrightarrow \limGs[(G,Y^*)]$ such that
\[
\alpha(\tau_x(t))=\tau_y(\alpha(t))\cdot s_x,\ \mbox{ for }\ t\in \limGs[(G,X^*)],
\]
where $s_x={h'}_y^{-1}\phi({g'}_y^{-1}g_x)h_x$ and $y$ is such that $g'_yg_x^{-1}\in H$.
\end{theor}
\begin{proof}
The statement follows from Sections~2.1--2.5 in \cite{self_sim_groups}.
\end{proof}

\textbf{Limit dynamical system.} The factor of the limit $G$-space $\limGs$ by the action of the
group $G$ is called the \textit{limit space} $\lims=\lims(\phi)$. It follows from the definition
that we may also consider $\lims$ as a factor of $\xmo$ by the following equivalence relation: two
left-infinite sequences $\ldots x_2 x_1,\ldots y_2 y_1$ are equivalent if and only if there exists
a left-infinite path $\ldots e_2 e_1$ in the nucleus $\nucl$ such that the edge $e_i$ is labeled by
$(x_i,y_i)$. The limit space $\lims$ is compact, metrizable, finite-dimensional space. It is
connected if the group $G$ is finitely generated and acts transitively on $X^n$ for all $n$.

The equivalence relation on $\xmo$ is invariant under the shift $\sigma$, therefore $\sigma$
induces a continuous surjective map $\si:\lims\rightarrow\lims$, and every point of $\lims$ has at
most $|X|$ preimages under $\si$. The dynamical system $(\lims,\si)$ is called the \textit{limit
dynamical system} of the self-similar action.

The image of $\xmo v$ for $v\in X^{n}$ in $\lims$ is called the \textit{tile} $\T_v$ of the $n$-th
level. Clearly
\[
\T_v=\bigcup_{x\in X}\T_{xv}\qquad\mbox{ and }\qquad \si(\T_{vx})=\T_v
\]
for every $v\in X^{*}$ and $x\in X$. Two tiles $\T_v$ and $\T_u$ of the same level have nonempty
intersection if and only if there exists $h\in\nucl$ such that $h(v)=u$. Under the open set
condition, every tile $\T_v$ is the closure of its interior, and any two different tiles of the
same level have disjoint interiors (Proposition~3.6.5 in \cite{self_sim_groups}). It will be used
in the next section that
\begin{equation}\label{eqn_diam_tile_goes_zero}
\lim_{n\rightarrow\infty} \max_{v\in X^n} \rm{diam} (\T_v)=0
\end{equation}
for any chosen metric on the limit space $\lims$ (see Theorem~3.6.9 in \cite{self_sim_groups}).

The inverse limit of the topological spaces $\lims\stackrel{\si}{\leftarrow} \lims
\stackrel{\si}{\leftarrow} \cdots$ is called the \textit{limit solenoid} $\sol$. One can consider
$\sol$ as a factor of the space $\xz$ of two-sided infinite sequences by the equivalence relation,
where two sequences $\xi,\eta$ are equivalent if and only if there exist a two-sided infinite path
in the nucleus labeled by the pair $(\xi,\eta)$. The two-sided shift on $\xz$ induces a
homeomorphism $\e:\sol\rightarrow\sol$.

Consider a simple example to illustrate the main notions of this section (see
\cite[Chapter~6]{self_sim_groups} for many other examples).

\textbf{Example.} Consider the transformation $a$ of the space $X^{*}$ over the binary alphabet
$X=\{0,1\}$ defined recursively by the rule
\begin{eqnarray*}
a(0w)=1w \quad \mbox{ and }\quad a(1w)=0a(w)
\end{eqnarray*}
for $w\in X^{*}$. The action of $a$ on the space $X^{*}$ corresponds to addition of one in the
binary numeration system, and therefore the transformation $a$ is called the (binary)
\textit{adding machine}. Indeed, $a(x_1x_2\ldots x_n)=y_1y_2\ldots y_n$ for $x_i,y_i\in X$ if and
only if
\[
1+x_1+2x_2+\ldots +2^{n-1}x_n=y_1+2y_2+\ldots+2^{n-1}y_n \mod 2^{n}.
\]
\begin{figure}
\begin{center}
\psfrag{a}{$a$} \psfrag{a1}{$a^{-1}$} \psfrag{id}{$id$} \psfrag{00}{$(0,0)$} \psfrag{11}{$(1,1)$}
\psfrag{01}{$(0,1)$} \psfrag{10}{$(1,0)$} \epsfig{file=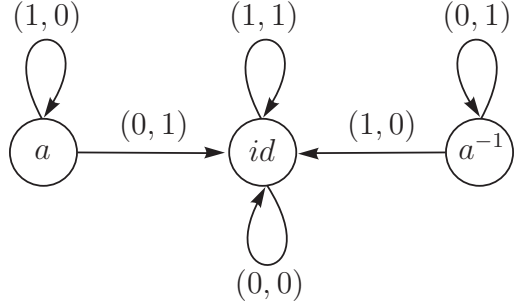,height=110pt}\caption{The nucleus
of the adding machine} \label{fig_AddingMachine}
\end{center}
\end{figure}
The group generated by $a$ is the infinite cyclic group, and we have a self-similar action of the
group $\Z\cong\langle a\rangle$ on the space $X^{*}$. This action corresponds to the virtual
endomorphism $\phi:2\Z\rightarrow\Z$ given by $\phi(m)=m/2$ and the coset transversal $T=\{0,1\}$
(there is no need to specify the cocycle for abelian groups). The nucleus of the action
$(\Z,X^{*})$ is shown in Figure~\ref{fig_AddingMachine}. It follows that two sequences in
$\xmo\times \Z$ are asymptotically equivalent if and only if they are of the form
\[
\ldots 000\cdot a^{n+1} \approx \ldots 111\cdot a^n \quad \mbox{or} \quad \ldots 001w\cdot a^n
\approx \ldots 110w\cdot a^n
\]
for $w\in X^{*}$ and $n\in\Z$. This is the usual identification of binary expansions of real
numbers (but written to the left). Hence the map $\Phi$ given by the rule
\[
\Phi(\ldots x_2x_1\cdot a^n)=n+\sum_{i=1}^{\infty} \frac{x_i}{2^i},\quad x_i\in X \mbox{ and }
n\in\Z,
\]
induces a homeomorphism between the limit $G$-space $\limGs[(\Z,X^{*})]$ and the real line
$\mathbb{R}$ with the standard action of $\Z$ by translations. The image of the tile $\tile$ under
$\Phi$ is the unit interval $[0,1]$. The limit space $\lims[(\Z,X^{*})]$ is homeomorphic to the
unit circle.



\section{Self-similar measures on limit spaces}\label{section_Measures}

Let us fix a contracting self-similar action $(G,X^*)$.\\

\noindent\textbf{Invariant measure on the limit $G$-space $\limGs$}. We consider the uniform
Bernoulli measure $\mu_u$ on the space $\xmo$ and the counting measure on the group $G$, and we put
the product measure on the space $\xmo\times G$. The push-forward of this measure under the factor
map $\pi_{\limGs[]}:\xmo\times G\rightarrow\limGs$ defines the \textit{measure $\mu$ on the limit
$G$-space} $\limGs$. The measure $\mu$ is a $G$-invariant $\sigma$-finite regular Borel measure on
$\limGs$.

\begin{prop}\label{prop_measure lims propert}
The measures of tiles have the following properties.
\begin{enumerate}
  \item $\mu(\tile_v)= |X|^n\cdot\mu(\tile_{uv})$ for every $v\in X^{*}$ and $u\in X^n$.
  \item $\mu(\tile_v\cap\tile_{v'})=0$ for $v,v'\in X^n, v\neq v'$.
  \item Let $\mu|_{\tile}$ be the measure $\mu$ restricted to the tile $\tile$. Then
\[
  \mu|_{\tile}(A)=\sum_{x\in X} \frac{1}{|X|}\mu|_{\tile}(\tau_x^{-1}(A))
\]
for any Borel set $A\subset\tile$.
\end{enumerate}
\end{prop}
\begin{proof}
Let us show that $\mu(A)\geq |X|\mu(\tau_x(A))$ for any Borel set $A$. Consider sequences which
represent points of the sets $A$ and $\tau_x(A)$:
\[
\pi^{-1}_{\limGs[]}(A)=\bigcup_{g\in G}T_g\cdot g \ \ \Rightarrow \ \
\pi^{-1}_{\limGs[]}(\tau_x(A))=\bigcup_{g\in G}T_g g(x)\cdot g|_x.
\]
It implies
\[
\mu(\tau_x(A))\leq \sum_{g\in G} \frac{1}{|X|}\mu_u (T_g)=\frac{1}{|X|}\mu(A).
\]

By applying this inequality $n$ times we get
$|X|^n\mu(\tile_v)=|X|^n\mu(\tau_v(\tile))\leq\mu(\tile)$ for $v\in{X}^n$. Since $\tile=\cup_{u\in
X^n}\tile_{u}$ we have that
\[
\mu(\tile)\leq\sum_{u\in X^n}\mu(\tile_{u})\leq\sum_{u\in X^n}\frac{1}{|X|^n}\mu(\tile)=\mu(\tile).
\]
Hence all the above inequalities are actually equalities, $\mu(\tile)=|X|^n\mu(\tile_{u})$ for
every $u\in X^n$, and $\mu(\tile_u\cap\tile_{u'})=0$ for different $u,u'\in X^n$.

Notice that since every Borel set $A\subset\tile$ can be approximated by unions of tiles of the
same level 
and using items 1 and 2 we have that if $\mu(\tau_x(A))<\varepsilon$ then $\mu(A)<\varepsilon|X|$.

It is left to prove item 3. First, let us show that the assertion holds for the tiles $\tile_v$.
Since the measure $\mu|_{\tile}$ is concentrated on the tile $\tile$, up to sets of measure zero
the set $\tau_x^{-1}(\tile_v)$ is equal $\tile_u$ if $v=ux$, and is empty if the last letter of $v$
is not $x$. Really, if $t\in\tau_x^{-1}(\tile_{ux})$ and $t\in\tile_v$ with $v\neq u, |v|=|u|$,
then $\tau_x(t)\in\tile_{vx}\cap\tile_{ux}$ and the measure of such points is zero. Hence
\[
\mu|_{\tile}(\tile_v)=\frac{1}{|X|}\mu|_{\tile}(\tile_u)=\sum_{x\in X}
\frac{1}{|X|}\mu|_{\tile}(\tau_x^{-1}(\tile_v)).
\]
Now we can approximate any Borel set by unions of tiles and pass to the limit.
\end{proof}

%
%
%

The tile $\tile$ of the limit $G$-space $\limGs$ can be considered as the attractor of the iterated
function system $\tau_x$, $x\in X$, i.e.
\[
\tile=\bigcup_{x\in X}\tau_x(\tile).
\]
Hence Proposition~\ref{prop_measure lims propert} item~3 implies that $\mu|_{\tile}$ is the
self-similar measure on $\tile$ by the standard definition of Hutchinson
(\ref{eqn_self_similar_measure}). The measure $\mu$ is the $G$-invariant extension of the
self-similar measure $\mu|_{\tile}$ to the limit $G$-space $\limGs$.


Let us show how to find the measure of the tile $\tile$. Let $\nucl$ be the nucleus of the action
$(G,X^*)$ identified with its Moore diagram. Replacing each label $(x,y)$ by label $x$ in the
nucleus $\nucl$ we get a strictly right-resolving graph denoted $\gr_{\nucl}$ labeled by elements
of $X$, so that we can apply the methods developed in Section~2.

\begin{theor}\label{theor_tile_measure}
The measure $\mu(\tile)$ is equal to the measure number $\mgr(\gr_\nucl)$ of the nucleus, in
particular it is always integer. Moreover, $\mu(\tile)=1$ if and only if the action satisfies the
open set condition.
\end{theor}
\begin{proof}
By Proposition~\ref{prop_asympt_Moore_diagr} we have
\[
\pi^{-1}_{\limGs[]}(\tile)=\bigcup_{g\in\nucl}F_g\cdot g^{-1},
\]
where the sets $F_g$ are defined using the graph $\gr_\nucl$ (see Section~2). Thus
\[
\mu(\tile)=\sum_{g\in\nucl}\mu_u(F_g)=\mgr(\nucl),
\]
which is integer by Theorem~\ref{theor_measure_graph_integer}. Observe that $\mu(\tile)\geq 1$
because $F_g=\xmo$ for $g=1\in\nucl$.

If the action satisfies the open set condition then $\mu(\tile)=\mgr(\gr_\nucl)=1$ by
Corollary~\ref{cor_meas_sum_meas}.

Suppose now that the action does not satisfy the open set condition. Then there exists an element
$h$ of the nucleus, whose all restrictions are non-trivial. Let $\nucl_1$ be the set of all
restrictions of $h$, and $\gr_{\nucl_1}$ be the corresponding graph. Then by Proposition~\ref{prop
cup B_v eq cup F_v}
\[
\sum_{g\in\nucl_1}\mu_u(F_g)=\mgr(\gr_{\nucl_1})\geq 1.
\]
Thus $\mu(\tile)=\sum_{g\in\nucl}\mu_u(F_g)\geq \mgr(\gr_{\nucl_1})+\mu_u(F_{g=1})\geq 2$.
\end{proof}

\begin{remark}
The measure $\mu(\tile)=\mgr(\gr_{\nucl})$ can be found algorithmically using the remarks after
Proposition~\ref{prop_eigenvect_transition_matr}.
\end{remark}


%
%

The next proposition shows that the covering $\limGs=\cup_{g\in G}\tile\cdot g$ is a perfect
multiple covering of multiplicity $\mu(\tile)$.

\begin{prop}\label{points_above}
Every point $x\in\limGs$ is covered by at least $\mu(\tile)$ tiles. The set $\dot{\limGs}$ of all
points $x\in\limGs$, which are covered by exactly $\mu(\tile)$ tiles, is open and dense in
$\limGs$, and its complement has measure $0$.
\end{prop}
\begin{proof}
For each $x\in\limGs$ we define the number $n_x$ of such $g\in G$ that the tile $\tile\cdot g$
contains $x$.

First we prove the inequality. Let $x\in\limGs$ be represented by the pair $w\cdot g$ from
$\xmo\times G$. Then $|V(w)|\geq\mgr(\gr_\nucl)$ by Theorem~\ref{theor_measure_graph_integer}. If
$h\in V(w)$ it means that the sequence $w$ ends in $h$, which by
Proposition~\ref{prop_asympt_Moore_diagr} means that there is a sequence $u_h\in\xmo$ such that $w$
is asymptotically equivalent to $u_h\cdot h$. It follows that $w\cdot g$ is asymptotically
equivalent to $u_h\cdot hg$. It means, in turn, that $x$ belongs to the tile $\tile\cdot hg$ for
every $h\in V(w)$. It follows that $n_x\geq\mgr(\gr_\nucl)$.

Consider the set $\mathcal{O}=\mathcal{O}(\gr_\nucl)$ defined in Lemma~\ref{lemma_set V open dense}
using the graph $\gr_\nucl$, in other words $\mathcal{O}$ is the set of all $w\in\xmo$ that end in
precisely $\mgr(\gr_\nucl)$ elements of $\nucl$. By the same considerations as above we see that if
$w\cdot g$ represents a point $x$ with $n_x=\mgr(\gr_\nucl)$ then $|V(w)|=\mgr(\gr_\nucl)$, that
is, $w$ belongs to $\mathcal{O}$. In other words, the set $\mathcal{O}\times G$ is closed under the
asymptotic equivalence relation on $\xmo\times G$, and it is the inverse image of the set
$\dot{\limGs}$ under the factor map $\pi_{\limGs[]}$. Since the set $\mathcal{O}\times G$ is open
and dense in $\xmo\times G$ by Lemma~\ref{lemma_set V open dense}, the same hold for $\dot{\limGs}$
in $\limGs$. The complement of $\dot{\limGs}$ has measure $0$, because the complement of the set
$\mathcal{O}$ has measure $0$ by Lemma~\ref{lemma_set V open dense}.
\end{proof}

\begin{cor}
Let $\partial\tile$ be the boundary of $\tile$. Then $\mu(\partial\tile)=0$.
\end{cor}
\begin{proof}
Take $x\in\partial\tile$. Suppose that $x$ is covered by exactly $k=\mu(\tile)$ tiles. Let
$U\subset\dot{\limGs}$ be the neighborhood of $x$ which intersects with finitely many tiles, say
with $l\geq k$ tiles. By subtracting any tile that does not contain $x$ we may suppose that $l=k$,
that is, $U$ has nonempty intersection only with those tiles that contain $x$. Now, since $x$ is in
the boundary of $\tile$, there is a point $x'$ in $U$ that is not covered by $\tile$ (otherwise
$x\in U\subset\tile$), thus $x'$ is covered by less than $k$ tiles, contradiction. Hence
$\partial\tile$ is in the complement of $\dot{\limGs}$ and has measure $0$.
\end{proof}


\begin{remark}\label{rem_meas_inters_tiles}
Two tiles $\tile\cdot g_1$ and $\tile\cdot g_2$ for $g_1,g_2\in G$ have nonempty intersection if
and only if $g_1g_2^{-1}\in\nucl$. Let us show how to find the measure of this intersection. By
Proposition~\ref{prop_asympt_Moore_diagr} we have
\begin{eqnarray*}
\pi^{-1}_{\limGs[]}((\tile\cdot g_1)\cap (\tile\cdot g_2))=\bigcup_{{\scriptsize\begin{array}{c}
                              h_1,h_2\in\nucl \\
                              g_1g_2^{-1}=h_1h_2^{-1}
                            \end{array}}}
(F_{h_1}\cap F_{h_2})\cdot h_1^{-1}g_1,
\end{eqnarray*}
where the sets $F_g$ are defined using the nucleus $\nucl$. The word problem in contracting
self-similar groups is solvable in polynomial time \cite[Proposition~2.13.10]{self_sim_groups} (one
can use the program package \cite{AutomGrp}). The measures of intersections $F_{h_1}\cap F_{h_2}$
can be found for example using method described after Corollary~\ref{cor_meas_gr_rational}.
\end{remark}


The measure $\mu=\mu_{(G,X^{*})}$ on the limit $G$-space $\limGs$ was defined using the specific
self-similar action $(G,X^{*})$ of the group $G$. Let us show that actually this measure depends
only on the associated virtual endomorphism $\phi$, as the limit $G$-space itself. It allows us to
consider the measure space $(\limGs(\phi),\mu)$ independently of the self-similar action. At the
same time, the measure $\mu(\tile)$ may vary for different self-similar actions as the nucleus
does. It is an interesting open question in what cases we can always choose a self-similar action
which satisfies the open set condition (see \cite{mult_an_self_sim,lag_wang:haar} for the abelian
case and applications to wavelets).

\begin{theor}\label{theor_measure_limGs_invariant}
Let $\phi: H\rightarrow G$ be a virtual endomorphism of the group $G$. Let $(G,X^*)$ and $(G,Y^*)$
be the contracting self-similar actions defined respectively by the triples $(\phi,T,C)$ and
$(\phi,T',C')$. Then the homeomorphism $\alpha : \limGs[(G,X^*)]\longrightarrow \limGs[(G,Y^*)]$
from Theorem~\ref{theor_limsX hom limsY} preserves measure, i.e.
\[
\mu_{(G,Y^*)}(\alpha(A))=\mu_{(G,X^*)}(A)
\]
for any Borel set $A$.
\end{theor}
\begin{proof}
Let $\nucl$ be the nucleus and $\tile$ be the tile of the action $(G,X^{*})$. By
Theorem~\ref{theor_tile_measure} and Theorem~\ref{theor_measure_graph_integer}
\begin{equation}\label{eqn_theor_meraXY_mu tile }
\mu_{(G,X^*)}(\tile)=\mgr(\gr_\nucl)=\min_{w\in\xmo} |V(w)|=k\in\mathbb{N},
\end{equation}
where $V(w)$ is defined using the graph $\gr_\nucl$. Consider
\[
\pi_Y^{-1}(\alpha(\tile))=\bigcup_{g\in G} T_g\cdot g\quad (\mbox{here } T_g\subseteq Y^{-\omega}),
\]
where $\pi_Y:Y^{-\omega}\times G\rightarrow \limGs$ is the canonical projection.

Take $w\in Y^{-\omega}$ and let us prove that there exist at least $k$ elements $g\in G$ such that
$w\in T_g$. The tiles $\tile\cdot g$ cover the limit $G$-space $\limGs[(G,X^*)]$, and we can find
$g\in G$ such that if $x=\alpha^{-1}(y)$ for $y=\pi_{Y}(w\cdot g)\in\limGs[(G,Y^*)]$ then $x$
belongs to the tile $\tile$. Then $x$ is represented by the sequence $u\cdot 1$ for $u\in\xmo$.
Equation~(\ref{eqn_theor_meraXY_mu tile }) implies $|V(u)|\geq k$ and so there exist $k$ elements
$h_1,\ldots, h_k\in\nucl$ and sequences $u_1,\ldots,u_k\in\xmo$ such that for every $i$ there
exists a left-infinite path in the nucleus $\nucl$ which ends in $h_i$ and is labeled by $(u,u_i)$.
By Proposition~\ref{prop_asympt_Moore_diagr}
\[
x=\pi_X(u\cdot 1)=\pi_X(u_i\cdot h_i)\quad\mbox{ and }\quad x\cdot h_i^{-1}\in\tile
\]
for every $i=1,\ldots, k$. Then
\[
\pi_Y(w\cdot gh_i^{-1})=\pi_Y(w\cdot g)h_i^{-1}=y\cdot h_i^{-1}=\alpha(x\cdot
h_i^{-1})\in\alpha(\tile)
\]
and thus $w\in T_{gh_i^{-1}}$ for all $i=1,\ldots, k$.

Let $\chi_{T_g}$ be the characteristic function of the set $T_g$. Then $\sum_{g\in
G}\chi_{T_g}(x)\geq k$ for almost all $x$. Integrating we get
$\mu_{(G,Y^*)}(\alpha(\tile))=\sum_{g\in G} \mu_u(T_g)\geq k=\mu_{(G,X^*)}(\tile)$.

Let us now show that $\mu_{(G,Y^*)}(\alpha(\tile_v))\leq\mu_{(G,X^*)}(\tile_v)$ for any $v\in X^*$.
Indeed, let $n=|v|$, then it follows from Theorem \ref{theor_limsX hom limsY} that
$\alpha(\tile_v)=\alpha(\tau_v(\tile))=\tau_u(\alpha(\tile))g$ for some word $u\in Y^n$ and $g\in
G$. Thus
\begin{eqnarray*}
\mu_{(G,Y^*)}(\alpha(\tile_v))&=&\mu_{(G,Y^*)}(\tau_u(\alpha(\tile)))\leq\frac{1}{|X|^n}\mu_{(G,Y^*)}(\alpha(\tile))=\\
&=&\frac{1}{|X|^n}\mu_{(G,X^*)}(\tile)=\mu_{(G,X^*)}(\tile_v).
\end{eqnarray*}

Let us prove that $\mu_{(G,Y^*)}\circ\alpha$ is absolutely continuous with respect to
$\mu_{(G,X^*)}$. Indeed, let $\mu_{(G,X^*)}(A)<\varepsilon$, $\pi^{-1}_{X}(A)=\cup_g T_g\cdot g$.
Then $\sum_g\mu_u(T_g)<\varepsilon$. It follows that there exist $v_{i,g}\in X^*$ such that
$T_g\subset\cup_i\xmo v_{i,g}$ and $\sum_{i,g}|X|^{-|v_{i,g}|}<\varepsilon$. Then
$A\subset\cup_{i,g}\tile_{v_{i,g}}\cdot g$, and we have that
$\sum_{i,g}\mu_{(G,X^*)}(\tile_{v_{i,g}})<\varepsilon\mu_{(G,X^*)}(\tile)$. Then
$\mu_{(G,Y^*)}(\alpha(A))\leq\sum\mu_{(G,Y^*)}(\alpha(\tile_{v_{i,g}}))\leq\sum\mu_{(G,X^*)}(\tile_{v_{i,g}})<\varepsilon\mu_{(G,X^*)}(\tile)$.
Since $\varepsilon$ is arbitrary, we are done.

We will now prove that $\mu_{(G,Y^*)}\circ\alpha\leq\mu_{(G,X^*)}$. Since both
$\mu_{(G,Y^*)}\circ\alpha$ and $\mu_{(G,X^*)}$ are invariant under multiplication by $g\in G$ it
suffices to prove this inequality for sets $A\subset\tile$. Since any Borel $A$ is a union of a
closed set and a set of arbitrarily small measure, it suffices to prove the inequality for closed
sets, as $\mu_{(G,Y^*)}\circ\alpha$ is absolutely continuous with respect to $\mu_{(G,X^*)}$ by
above.

So let $A\subset\tile$ be a closed set. For each $n$, let $A_n$ be the union of all tiles
$\tile_v,v\in X^n,$ that have non-empty intersection with $A$. Let us show that $A=\cap_n A_n$.
Suppose $x\in\cap A_n$, and $x$ is not in $A$. Then for each $n$ there is $v_n\in X^n$ such that
$\tile_{v_n}$ has nonempty intersection with $A$ and $x\in\tile_{v_n}$. It follows that $x$ has
some representation $u_nv_n\in\xmo$. Since the number of such representations is finite, we may
choose subsequence $n_k$ such that $v_{n_k}$ is the beginning of some word $v\in\xmo$. Since $A$ is
compact it follows that $A$ has nonempty intersection with $\cap_k\tile_{v_{n_k}}$, thus
$\cap_k\tile_{v_{n_k}}$ contains at least two points, which is impossible.

Since $A_n$ is the union of some tiles of level $n$, which are disjoint up to sets of measure $0$,
the inequality also holds for all $A_n$. Going to the limit, we get the inequality for $A$.

By interchanging $X$ and $Y$ we get the reverse inequality, and we are done.
\end{proof}

\begin{remark}
It is important in the theorem that we take the uniform Bernoulli measure on $\xmo$. The problem is
that the homeomorphism $\alpha$ may change the Bernoulli measure with a non-uniform weight to a
measure that is not Bernoulli.
\end{remark}

\vspace{0.3cm}\noindent\textbf{Self-similar measure on the limit space $\lims$}. The push-forward
of the uniform Bernoulli measure $\mu_u$ under the factor map $\pi_{\lims[]}:\xmo\rightarrow\lims$
defines the \textit{self-similar measure $\muls$ on the limit space} $\lims$. The measure $\muls$
is a regular Borel probability measure on $\lims$. The shift $\si$ is a measure-preserving
transformation of $\lims$.

Consider the set $\mathcal{U}=\mathcal{U}(\nucl)$ of all sequences $w\in\xmo$ with the property
that every left-infinite path in the nucleus $\nucl$ labeled by $(w,w)$ ends in $1$. Define
$\mathcal{U}^*$ as the set of finite words that satisfy the same condition.

\begin{lemma}\label{lemma U_open_dense_meas1}
The set $\mathcal{U}$ is open and dense in $\xmo$, and $\mu_p(\mathcal{U})=1$. For each $w\in
\mathcal{U}$ there is a beginning of $w$ that belongs to $\mathcal{U}^*$, and
$\mathcal{U}=\cup_{w\in\mathcal{U}^*}\xmo w$.

The sets $\mathcal{U}$ and $\mathcal{U}\times G$ are closed under the asymptotic equivalence
relation on $\xmo$ and $\xmo\times G$ respectively.
\end{lemma}
\begin{proof}
Construct the graph $\gr$ with the elements of $\nucl$ as vertices, for each edge labeled by
$(x,x)$ in $\nucl$ we have the edge in $\gr$ with the same starting and end vertices, and labeled
by $x$; and for each edge in $\nucl$ labeled by $(x,y)$ for $x\neq y$ there is an edge in $\gr$
with the same starting vertex that ends in the trivial element and labeled by $x$. Then the set
$\mathcal{U}$ coincides with the set $\mathcal{O}(\gr)$. Indeed, $h\in V_{\gr}(w)$ if and only if
there is a path in $\nucl$ labeled by $(w,w)$ that ends in $h$, thus $w\in\mathcal{U}$ if and only
if $V_{\gr}(w)=\{1\}$. For every nontrivial element $h\in\nucl$ there exists a word $v\in X^{*}$
such that $h(v)\neq v$. It follows that there exists a path in the graph $\gr$ from $h$ to $1$.
Hence the component $\{1\}$ is the only strongly connected component of the graph $\gr$ without
outgoing edges. By Corollary~\ref{cor_meas_sum_meas} the measure number of $\gr$ is $1$, and
$\mathcal{U}=\mathcal{O}(\gr)$. The first statement of the lemma now follows from
Lemma~\ref{lemma_set V open dense}.

Let us show that the set $\mathcal{U}$ is closed under the asymptotic equivalence relation (then
the set $\mathcal{U}\times G$ is also closed). It is sufficient to show that if there is a path in
$\nucl$ labeled by $(u,v)$ and $u\in\mathcal{U}$ ($v\in\mathcal{U}$) then $v\in\mathcal{U}$
($u\in\mathcal{U}$). Let the path in $\nucl$ labeled by $(u,v)$ end in $h$. It follows that $u$ is
asymptotically equivalent to $v\cdot h$. Suppose there is a path in $\nucl$ labeled by $(v,v)$ that
ends in $g$. It follows that $v$ is asymptotically equivalent to $v\cdot g$. Thus $u$ is
asymptotically equivalent to $v\cdot gh$ which is asymptotically equivalent to $u\cdot h^{-1}gh$.
By definition, there is a path in $\nucl$ labeled by $(u,u)$ which ends in $h^{-1}gh$. Since
$u\in\mathcal{U}$ we get $h^{-1}gh=1$, thus $g=1$.
\end{proof}


\begin{prop}\label{prop_measure_prop_lims}
Almost every point of $\lims$ has precisely $|X|$ preimages under~$\si$.
\end{prop}
\begin{proof}
Since every point of $\lims$ has at most $|X|$ preimages under $\si$, it is enough to show that for
almost every $w\in\xmo$ the map $\pi_{\lims[]}:\sigma^{-1}(w)\rightarrow\si^{-1}(\pi_{\lims[]}(w))$
is one-to-one. Suppose that for some $w\in\xmo$ and $x\neq y$ in $X$ we have
$\pi_{\lims[]}(wx)=\pi_{\lims[]}(wy)$. It follows that $wx,wy$ are asymptotically equivalent, thus
there is a left-infinite path in $\nucl$ labeled by $(wx,wy)$. It follows that the prefix of this
path labeled by $(w,w)$ must end in the nontrivial element, so $w\not\in\mathcal{U}$. Since
$\mu_u(\mathcal{U})=1$ by Lemma~\ref{lemma U_open_dense_meas1}, we are done.
\end{proof}


The definition of the measure $\muls$ uses encoding of $\lims$ by sequences $\xmo$. At the same
time the space $\lims$ is defined as the space of obits $\limGs/G$. Let us show that we can recover
$(\lims,\muls)$ from the measure space $(\limGs,\mu)$.

\begin{prop}\label{x_j_homeo}
Let $\rho$ be the factor map $\limGs\rightarrow\lims$. Then for any $u\in\mathcal{U}^*$ the
restriction $\rho|_{\tile_u}:\tile_u\rightarrow\T_u$ is a homeomorphism,
$\rho^{-1}(\T_u)=\sqcup_{g\in G}\tile_u\cdot g$.
\end{prop}
\begin{proof}
By definition, the map $\rho|_{\tile_u}:\tile_u\rightarrow\T_u$ is surjective, and $\tile_u$ is
compact. Hence in order to prove that $\rho$ is a homeomorphism it is left to show that $\rho$ is
injective on $\tile_u$. Take $x,y\in\tile_u$, and let $wu\in\xmo$ represent $x$ and $vu\in\xmo$
represent $y$. Suppose that $\rho(x)=\rho(y)$. By Proposition~\ref{prop_asympt_Moore_diagr} it
means that there is a left-infinite path in $\nucl$ labeled by $(wu,vu)$. Since
$u\in\mathcal{U}^*$, this path must end in $1$. It follows, that $wu$ and $vu$ represent the same
point of the tile $\tile_u$, and $x=y$.

To prove the second claim, take $x\in\tile_u\cdot g\cap\tile_u\cdot g'$. Then $x$ is represented by
two asymptotically equivalent sequences $wu\cdot g$ and $w'u\cdot g'$. It follows that there is a
path in the nucleus $\nucl$ labeled by $(u,u)$ which ends in $g'g^{-1}$. Then $g=g'$ since
$u\in\mathcal{U}^*$.
\end{proof}

\begin{theor}
The projection $\limGs\rightarrow\lims$ is a covering map up to sets of measure zero.
\end{theor}
\begin{proof}
Consider the sets $\tilde{\limGs}=\pi_{\limGs[]}(\mathcal{U}\times G)$ and
$\tilde{\lims}=\pi_{\lims[]}(\mathcal{U})$. It follows that $\tilde{\limGs}/G=\tilde{\lims}$. Since
the set $\mathcal{U}$ has measure $1$ by Lemma~\ref{lemma U_open_dense_meas1}, the complements of
$\mathcal{U}\times G$, of $\tilde{\limGs}$, and of $\tilde{\lims}$ have measure $0$.


Since the group $G$ acts properly on $\limGs$, the same holds for $\tilde{\limGs}$. It is left to
prove the freeness. Suppose $x\cdot g=x$ for $x\in\tilde{\limGs}$ and $g\in G$. Let $u\cdot h$ be a
representative of $x$ in $\xmo\times G$. It follows that $u\cdot hg$ is asymptotically equivalent
to $u\cdot h$, thus there is a path in $\nucl$ labeled by $(u,u)$ that ends in $h^{-1}gh$. Since
$u\in\mathcal{U}$ we have $g=1$.

Hence the projection $\tilde{\limGs}\rightarrow\tilde{\lims}$ is a covering map, and the statement
follows.
\end{proof}

Let $(\limGs[],\nu)$ be a locally compact measure space, the group $G$ acts freely and properly
discontinuously on $\limGs[]$ by homeomorphisms, and the measure $\nu$ is $G$-invariant. There is a
unique measure $\nu_{*}$ on the quotient space $\limGs[]/G$, called the \textit{quotient measure},
with the property that if $U$ is an open subset of $\limGs[]$ such that $Ug\cap Uh=\emptyset$ for
all $g,h\in G$, $g\neq h$, then $\nu_{*}(U/G)=\nu(U)$.


\begin{prop}
The quotient measure $\mu_{*}$ of the limit $G$-space $(\limGs,\mu)$ coincides with the measure
$\muls$ on the limit space $\lims$.
\end{prop}
\begin{proof}
Consider the sets $\tilde{\lims}$ and $\tilde{\limGs}$ of full measure from the previous theorem.
For every $u\in\mathcal{U}^*$ we have $\T_u\subset\tilde{\lims}$, $\tile_u\cdot
g\subset\tilde{\limGs}$ for every $g\in G$, and $\rho^{-1}(\T_u)=\bigsqcup_{g\in G} \tile_u\cdot g$
by Proposition~\ref{x_j_homeo}. Since $\tilde{\lims}=\cup_{u\in\mathcal{U}^*}\T_u$ by
Lemma~\ref{lemma U_open_dense_meas1} it suffices to show that for any $u\in\mathcal{U}^*$ the
restriction of $\mu_{*}$ on $\T_u$ is equal to the restriction of $\muls$. Take a Borel set
$A\subset\tile_u$ and consider its preimage $\pi^{-1}_{\limGs[]}(A)=\sqcup_{g\in\nucl} A_g\cdot g$.
Then $\pi^{-1}_{\lims[]}(\rho(A))=\sqcup_{g\in\nucl} A_g$. Here the union is disjoint because if
$w\in{A_g\cap{}A_h}$ then $w\cdot g$ and $w\cdot h$ are asymptotically equivalent to $vu$ and $v'u$
respectively, for some $v,v'\in\xmo$. It follows that $vu\cdot g^{-1}$ is equivalent to $v'u\cdot
h^{-1}$. Since $u\in\mathcal{U}$, $g=h$. Hence
\[
\muls(\rho(A))=\mu_u(\pi^{-1}_{\lims[]}(\rho(A)))=\mu_u(\sqcup_{g\in\nucl}
A_g)=\sum_{g\in\nucl}\mu_u(A_g).
\]
By the property of the quotient measure
\[
\mu_{*}(\rho(A))=\mu(A)=\sum_{g\in\nucl}\mu_u(A_g)=\muls(\rho(A)).
\]
\end{proof}

\begin{cor}
Theorem~\ref{theor_measure_limGs_invariant} holds for the limit space $\lims$, i.e. the measure
space $(\lims,\muls)$ depends only on the associated virtual endomorphism.
\end{cor}

\begin{cor}
$\muls(\T_u)=\mu(\tile_u)=\frac{1}{|X|^{|u|}}\mgr(\tile)$ for $u\in\mathcal{U}^*$.
\end{cor}

%
%
%

\begin{theor}\label{theor_lims_conju_Bernoulli}
$(\lims,\si,\muls)$ is conjugate to the one-sided Bernoulli $|X|$-shift.
\end{theor}
\begin{proof}
We will use notions and results from \cite{hoffman:Bernoulli}. First recall that a measure
preserving map with entropy $\log{d}$ is called {\it uniformly $d$-to-one endomorphism} if it is
almost everywhere $d$-to-one and the conditional expectation of each preimage is $1/d$.  The
standard example is $(\xmo,\sigma,\mu_u)$, which is uniformly $|X|$-to-one. Next we want to show
that $(\lims,\si,\muls)$ is also uniformly $|X|$-to-one. By
Proposition~\ref{prop_measure_prop_lims} the map $\pi:\xmo\rightarrow\lims$ is injective on the
preimages $\sigma^{-1}(w)$ for almost all $w\in\xmo$, that is, $\pi$ is {\it tree adapted} in
terminology of \cite{hoffman:Bernoulli}. We can apply Lemma~2.3 from \cite{hoffman:Bernoulli},
which says that a tree adapted factor of a uniform $d$-to-one endomorphism is again uniform
$d$-to-one endomorphism. In particular, the shift $\si$ is the map of maximal entropy $\log |X|$.

To prove the theorem, we use the following Theorem~5.5 in \cite{hoffman:Bernoulli}.

\begin{theor}
A uniform $d$-to-one endomorphism $(Y,S,\mu)$ is one-sidedly conjugated to the one-sided Bernoulli
$d$-shift if and only if there exists a generating function $f$ so that $(Y,S,\mu)$ and $f$ are
tree very weak Bernoulli.
\end{theor}


Recall the definition of tree very weak Bernoulli and generating function. Let $(Y,S,\mu)$ be
uniformly $d$-to-one and $f:Y\rightarrow{R}$ be a tree adapted function to a compact metric space
$R$ with metric $D$. The function $f$ is called \textit{generating} if the $\sigma$-algebra on $Y$
is generated by $S^{-i}f^{-1}(\mathcal{B}),i\geq{0}$, where $\mathcal{B}$ is the $\sigma$-algebra
of Borel sets of the space $R$.

Informally, ``tree very weak Bernoulli'' means that for almost all pairs of points in $Y$ their
trees of preimages are close. To give a formal definition note that since $S$ is uniformly
$d$-to-one, for almost all points $y\in Y$ the set $S^{-k}(y)$ contains exactly $d^k$ points, i.e.
the tree of preimages is a $d$-regular rooted tree. The set $\{1,\ldots, d\}^{*}$ of finite words
over $\{1,\ldots,d\}$ can be considered as a $d$-regular rooted tree, where every word $v$ is
connected with $vx$ for $x\in\{1,\ldots,d\}$, and the root is the empty word $\emptyset$. We can
use the tree $\{1,\ldots, d\}^{*}$ to label the trees of preimages. For almost all points $y\in Y$
there is a map $T_y:\{1,\dots,d\}^*\rightarrow{Y}$ such that $T_y(\emptyset)=y$ and
$T_y(\sigma(v))=S(T_y(v))$ for all nonempty words $v\in \{1,\dots,d\}^*$. Every map $T_y$ is tree
adapted, and it is uniquely defined up to an automorphism of the tree $\{1,\dots,d\}^*$. Then
$(Y,S,\mu)$ and $f$ are called {\it tree very weak Bernoulli} if for any $\varepsilon>0$ and all
sufficiently large $n$ there is a set $W=W(\varepsilon,n)\subset Y$ with $\mu(W)>1-\varepsilon$
such that for any $y,y'\in{W}$
\begin{equation}
\label{} t_n(y,y')=\min_\psi \frac{1}{n}\sum_{v\in \{1,\dots,d\}^{*}, |v|\leq
n}d^{-|v|}D(f(T_y(v)),f(T_{y'}(\psi(v))))<\varepsilon,
\end{equation}
where the minimum is taken over all automorphisms $\psi$ of the tree $\{1,\dots,d\}^*$. Notice that
the definition of $t_n$ does not depend on the choice of $T_y$.

Let us show that $(\lims,\si,\muls)$ is tree very weak Bernoulli for the identity map
$id:\lims\rightarrow\lims$. It immediately follows that $id$ is tree adapted and generating. Take a
point $x\in\lims$ and let $x$ be represented by some $w\in\xmo$. Define the map
$T_x:X^{*}\rightarrow\lims$ by the rule $T_x(v)=\pi(wv)$. It is enough to show that $
t_n(x,x')\rightarrow 0, n\rightarrow\infty, $ for almost all $x,x'\in\lims$. Using
(\ref{eqn_diam_tile_goes_zero}) we can find $n_1$ such that $\max_{v\in
X^n}\rm{diam}(\T_v)<\epsilon/2$ for all $n\geq n_1$. It means that $D(\pi(wv),\pi(w'v))<\epsilon/2$
for all $v,|v|\geq n_1$ (here $D$ is a fixed metric on the limit space $\lims$). Thus, taking
$\psi$ to be the identical tree automorphism, we have that
\[
t_n(x,x')<\frac{n_1}{n} \rm{diam}(\lims)+\frac{n-n_1}{n}\frac{\epsilon}{2}<\frac{n_1}{n}
\rm{diam}(\lims)+\frac{\epsilon}{2}<\epsilon
\]
for $n>2n_1\rm{diam}(\lims)/\epsilon$. Hence $(\lims,\si,\muls)$ and $id$ are tree very weak
Bernoulli, which finishes the proof.
\end{proof}



In the same way we introduce the measure $\muls_\e$ on the limit solenoid $\sol$ as the
push-forward of the uniform Bernoulli measure on $\xz$. It is easy to see that $(\sol,\e,\muls_\e)$
is the inverse limit of dynamical systems $(\lims,\si,\muls)$ (see \cite[page~27]{dyn_sysII}). In
particular, we get

\begin{cor}
$(\sol,\e,\muls_\e)$ is conjugate to the two-sided Bernoulli $|X|$-shift.
\end{cor}

\section{Applications and Examples}\label{section_Applictions}

\textbf{Invariant measures on nilpotent Lie groups.} Let $G$ be a finitely generated torsion-free
nilpotent group. Let $\phi:H\rightarrow G$ be a contracting surjective virtual endomorphism such
that the associated self-similar action is faithful (i.e. $\phi$-core($H$) is trivial in the
terminology of \cite{sidki:virt_nilp}). Then $\phi$ is also injective by Theorem~1 in
\cite{sidki:virt_nilp}, thus $\phi$ is an isomorphism and we can apply Theorem~6.1.6 from
\cite{self_sim_groups}. The group $G$ and its subgroup $H$ are uniform lattices of a simply
connected nilpotent Lie group $L$ by Malcev's completion theorem. The isomorphism
$\phi:H\rightarrow G$ extends to a contracting automorphism $\phi_L$ of the Lie group $L$. There
exists a $G$-equivariant homeomorphism $\Phi:\limGs\rightarrow L$ such that
$\phi_L(\Phi(t))=\Phi(\tau_{x_0}(t)\cdot g_0)$ for every $t\in\limGs$ and fixed $x_0\in X$ and
$g_0\in G$.





\begin{prop}\label{haar}
The push-forward $\Phi_*\mu$ of the measure $\mu$ on the limit $G$-space $\limGs$ is the (right)
Haar measure on the Lie group $L$.
%
\end{prop}
\begin{proof}
The measure $\Phi_*\mu$ is a non-zero regular Borel measure on $L$. It is left to prove that it is
translation invariant. Since the measure $\mu$ is $G$-invariant and the map $\Phi$ is
$G$-equivariant, the measure $\Phi_*\mu$ is $G$-invariant. By the property of the map $\Phi$ we
have $\phi_L(\Phi(A))=\Phi(\tau_{x_0}(A))g_0$ for every Borel set $A\subset\limGs$. Notice that
since the map $\phi$ is injective, we get $\mu(A)=|X|\mu(\tau_x(A))$ (see the proof of
Proposition~\ref{prop_measure lims propert}) and hence $\Phi_*\mu(B)=|X|\Phi_*\mu(\phi_L(B))$ for
every Borel set $B\subset L$. It follows that the measure $\Phi_*\mu$ is $\cup_n
\phi_L^n(G)$-invariant. Since $\phi_L$ is contracting, the set $\cup_n \phi_L^n(G)$ is dense in the
Lie group $L$. Hence $\Phi_*\mu$ is $L$-invariant.
\end{proof}


The same observation holds in a more general setting of finitely generated virtually nilpotent
groups (or under the conditions of Theorem~6.1.6).

\vspace{0.3cm}\noindent\textbf{Lebesgue measure of self-affine tiles.} Let $A$ be an $n\times n$
integer expanding matrix, where expanding means that every eigenvalue has modulus $>1$. The lattice
$\Z^n$ is invariant under $A$, and we can choose a coset transversal $D=\{d_1,\ldots,d_m\}$ for
$\Z^n/A(\Z^n)$, where $m=|det(A)|$.
There exists a unique nonempty compact set $T=T(A,D)\subset\mathbb{R}^n$, called (standard)
\textit{integral self-affine tile}, satisfying
\[
A(T)=\bigcup_{d\in D} T+d.
\]
The tile $T$ has positive Lebesgue measure, is the closure of its interior, and the union above is
nonoverlapping (the sets have disjoint interiors) \cite{lag_wang:self_affine}. It is well known
\cite{lag_wang:int_self_affine_I} that the tile $T$ has integer Lebesgue measure. The question how
to find this measure is studied in
\cite{lag_wang:int_self_affine_I,he_hui:self_affine,deng_he:int_self_affine}, and finally answered
in \cite{gab_yu:natur_til}. A related question how to find the measure of intersection $T\cap
(T+a)$ for $a\in\mathbb{Z}^n$ is studied in \cite{gab_yu:natur_til,int_two_tiles}. Let us show how
to answer these questions using the theory of self-similar groups.

The inverse of the matrix $A$ can be considered as the contracting virtual endomorphism $A^{-1}:
A(\Z^n)\rightarrow\Z^n$ of the group $\Z^n$, which is actually an isomorphism so that we can apply
the previous example of this section. Put $X=\{x_1,\ldots,x_m\}$ and let $(\mathbb{Z}^n,X^{*})$ be
the self-similar contracting action defined by the virtual endomorphism $A^{-1}$, the coset
transversal $D$, and the trivial cocycle $C=\{1,\ldots,1\}$ (see Section~3). The group $\Z^n$ is
the uniform lattice in the Lie group $\mathbb{R}^n$. Hence by Theorem~6.1.6 in
\cite{self_sim_groups} (see also Section~6.2 there) there exists a $\Z^n$-equivariant homeomorphism
$\Phi:\limGs[\Z^n]\rightarrow\mathbb{R}^n$ given by
\[
\Phi(\ldots x_{i_2} x_{i_1}\cdot g)=g+A^{-1}d_{i_1}+A^{-2}d_{i_2}+\ldots
\]
for $i_j\in\{1,\ldots,m\}$ and $g\in\Z^n$. The image of the tile $\tile$ is the self-affine tile
$T$.

%
%

\begin{prop}
The push-forward $\Phi_*\mu$ of the measure $\mu$ on the limit $G$-space $\limGs[\Z^n]$ is the
Lebesgue measure $\lambda$ on $\mathbb{R}^n$.
\end{prop}
\begin{proof}
The measure $\Phi_*\mu$ is the Haar measure on $\mathbb{R}^n$ by the above example. Since the Haar
measure is unique up to multiplicative constant, we have that $\Phi_*\mu=c\lambda$ for some
constant $c>0$ and we need to prove that $c=1$.

Recall that $\mu(\tile)$ is integer by Proposition~\ref{theor_tile_measure}, and almost every point
of $\limGs$ is covered by $\mu(\tile)$ tiles $\tile\cdot g$ by Proposition~\ref{points_above}.
Hence $\Phi_*\mu(T)=\Phi_*\mu(\Phi(\tile))=\mu(\tile)$ is integer and almost every point of
$\mathbb{R}^n$ is covered by $\mu(\tile)$ tiles $T+g$, $g\in\Z^n$, with respect to the measure
$\Phi_*\mu$, and thus with respect to the Lebesgue measure $\lambda$. It follows that, if
$\chi_{T+g}$ is the characteristic function of $T+g$, then $\sum_{g\in\Z^n}\chi_{T+g}=\mu(\tile)$
almost everywhere with respect to both measures. Hence
\[
\Phi_*\mu(T)=\mu(\tile)=\int_I\sum_{g\in\Z^n}\chi_{T+g}d\lambda=\sum_{g\in\Z^n}\int_{I+g}\chi_T
d\lambda=\int_{\mathbb{R}^n}\chi_T d\lambda=\lambda(T),
\]
where $I$ is the unit cube in $\mathbb{R}^n$. Since $\mu(\tile)$ is positive, $c=1$.
\end{proof}


\begin{cor}
The Lebesgue measure of the self-affine tile $T$ is equal to the measure number $\mgr(\nucl)$ of
the nucleus $\nucl$ of the associated self-similar action $(\mathbb{Z}^n,X^{*})$.
\end{cor}

The nucleus of a contracting self-similar action can be found algorithmically using the program
package \cite{AutomGrp}, and the number $\mgr(\nucl)$ can be found using the remarks after
Proposition~\ref{prop_eigenvect_transition_matr}. The measures of sets $T\cap (T+a)$ for
$a\in\mathbb{Z}^n$ can be found by Remark~\ref{rem_meas_inters_tiles}.



The methods developed in \cite{gab_yu:natur_til} to find the Lebesgue measure of integral
self-affine tiles are related to the discussion above. Take the complete automaton
\cite[page~11]{self_sim_groups} of the self-similar action $(\Z^n,X^{*})$ (it actually coincides
with the graph $B(\Z^n)$ from \cite{canon_numsys}), revert the direction of every edge, identify
edges with the same starting and end vertices labeled by $(d_1,d_2)$ with the same difference
$r=d_2-d_1$, and put the new label $r$ on this edge. We get the graph $\mathcal{G}(\Z^n)$
constructed in \cite{gab_yu:natur_til} and \cite{fract_numsys}. The set $W$ constructed using
$\mathcal{G}(\Z^n)$ \cite[page~195]{gab_yu:natur_til} is precisely the nucleus $\nucl$. Hence the
theory of self-similar groups provides a nice explanation to the ideas in
\cite[Section~3]{gab_yu:natur_til} and the methods developed in Sections~2,4 can be seen as its
non-abelian generalization.

It is shown in \cite{lag_wang:int_self_affine_II} that integral self-affine tile $T$ gives a
lattice tiling  of $\mathbb{R}^n$ with some lattice $L\subset\Z^n$. An interesting open question is
whether this holds for any (self-replicating) contracting self-similar action $(G,X^{*})$ (or at
least for self-similar actions of torsion-free nilpotent groups), i.e. the tile $\tile$ gives a
tiling of $\limGs$ with some subgroup $H< G$.

\vspace{0.3cm}\noindent\textbf{Invariant measures on Julia sets of rational functions.} Let
$f(z)\in\mathbb{C}(z)$ be a rational map of the Riemann sphere $\Rs$ of degree $d\geq 2$. It is
proved in \cite{invmeas_rat} that there exists a natural $f$-invariant probability measure $\mu_f$
with the support on the Julia set $J(f)$. In \cite{rational_bernoulli}, D.~Heicklen and C.~Hoffman
showed that the dynamical system $(J(f),f,\mu_f)$ is conjugated to the one-sided Bernoulli
$d$-shift, proving the conjecture in \cite{invmeas_rat,lub:entr_ration}. In the case of
sub-hyperbolic rational functions this result also follows from
Theorem~\ref{theor_lims_conju_Bernoulli} using the theory of iterated monodromy groups developed in
\cite{self_sim_groups}.

Let $f$ be a sub-hyperbolic rational function and let $\M=\Rs\setminus\overline{P}_f$, where $P_f$
is the post-critical set of $f$. Then $f$ defines a $d$-fold partial self-covering $f:
f^{-1}(\M)\rightarrow\M$. Take a base point $t\in\M$ and let $T_t$ be the tree of preimages
$f^{-n}(t)$, $n\geq 0$, where every vertex $z\in f^{-n}(t)$ is connected by an edge to $f(z)\in
f^{-n+1}(t)$. The fundamental group $\pi_1(\M,t)$ acts by automorphisms on $T_t$ through the
monodromy action on every level $f^{-n}(t)$. The quotient of $\pi_1(\M,t)$ by the kernel of its
action on $T_t$ is the \textit{iterated monodromy group} $IMG(f)$ of the map $f$. The group
$IMG(f)$ is contracting self-similar group and the limit dynamical system $(\lims[IMG(f)],\si)$ is
conjugated to $(J(f),f)$ by Theorem~6.4.4 in \cite{self_sim_groups}. Under this conjugation the
measure $\muls$ is send to the measure $\mu_f$, and we can apply
Theorem~\ref{theor_lims_conju_Bernoulli} to the dynamical system $(J(f),f,\mu_f)$.


\bibliographystyle{plain}

\end{document}